\documentclass[letter,11pt]{article}
\usepackage{amsfonts}
\usepackage{amssymb}
\usepackage{amsmath}
\usepackage{cite}

\newtheorem{prop}{Proposition}

\newtheorem{lema}[prop]{Lemma}

\newcommand{\R}{\mathbb{R}}
\newcommand{\D}{\mathbb{D}}
\newcommand{\T}{\mathbb{T}^1}

\newcommand{\C}{\mathbb{C}}

\newcommand{\Z}{\mathbb{Z}}

\newcommand{\cc}{\mathcal{C}}

\title{An example showing that fibred quadratic polynomials admit many attracting invariant curves}
\author{Mario PONCE\footnote{Partially supported by FONDECYT 3080055.}\\
PUC-Chile\\
mponcea@mat.puc.cl}
\begin{document}
\maketitle
\begin{abstract}
We present an example of a fibred quadratic polynomial admitting an attracting invariant  $2$-curve. By an unfolding construction we obtain an example of a fibred quadratic polynomial admitting two attracting invariant curves. This phenomena can not occur in the non-fibred setting. 
\end{abstract}
\section{Introduction}
A classical result due to Fatou and Julia asserts that quadratic polynomials may have at most one attracting cycle. This is an easy consequence of the fact that every immediate basin of attraction must contain a critical point (see \cite{CAGA93}). In this note we deal with a more general version of quadratic dynamics, namely,   fibred quadratic polynomials over an irrational rotation (also known as quasi periodically forced quadratic maps).  In \cite{PONC07} the author shows that many features of the classical theory of local dynamics also hold for the fibred setting. In \cite{SEST99} Sester defines  hyperbolicity for this class of maps and shows many equivalent characterizations which are analogous to the non-fibred (classical) characterizations of hyperbolicity for quadratic polynomials. 
\\

In this short note we show a new phenomena on fibred quadratic polynomials that do not occur in the non-fibred setting. We show an example of a fibred quadratic polynomial admitting two attracting invariant curves.  The construction follows a simple scheme: first, we present a fibred quadratic polynomial with an attracting invariant  $2$-curve. By an unfolding construction we obtain the desired example. 
\section{Definitions}

  A {\it fibred quadratic polynomial} is a continuous map  
  \begin{eqnarray*}
    P:\T\times \C&\to &\T\times \C \\
  (\theta,z)&\mapsto&\big(\theta+\alpha,P_{\theta}(z)\big)
  \end{eqnarray*}
  where $\alpha$ is a real angle (in $\T$) and $P_{\theta}:\C\to \C$ is a quadratic polynomial varying continuously with $\theta$. In general we will take $P_{\theta}(z)=P^{\cc}_{\theta}(z)=z^2+\cc(\theta)$, where $\cc:\T\to\C$ is a continuous function. The function $\cc$ must be thought as a parameter. This family  of fibred quadratic polynomials has been widely studied by Sester \cite{SEST99}. In that article, Sester defines the Julia and filled in Julia sets, studies the Green function and some kind of hyperbolicity on the Julia set. Sester also identifies in several equivalent ways the parameters $\cc$ which give raise to quasi-circle Julia set. That is, Sester defines and characterizes the {\it main cardioid} of the connectedness locus for the family $\{P^{\cc}\}_{\cc:\T\to\C}$.  \\
  
  Since $\alpha$ is irrational  the map $P$ has no fixed points (nor periodic orbits). A natural invariant object that extend the notion  of a fixed point is  an {\it invariant curve}, that is, a continuous curve $\gamma:\T\to \D$ such that $$P(\theta, (\gamma(\theta))=(\theta+\alpha, \gamma(\theta+\alpha)).$$ The real number 
  \[
  \kappa(\gamma)=\int_{\T}\log |\partial_zP_{\theta}\big(\gamma(\theta)\big)|d\theta
  \]
  is called the {\it multiplicator} (or Lyapunov exponent in the complex direction) and measures the expanding or contracting infinitesimal nature of the invariant curve. In \cite{PONC07} the author shows that an expanding or contracting multiplicator (that is $\kappa(\gamma)>1$ or $\kappa(\gamma)<1$) implies actually that the map is conjugated to a linear fibred map which is  expanding or contracting respectively.  
\\

Let  $\tilde \gamma:[0,2]\to \C$ be a continuous map such that $\tilde\gamma(0)=\tilde\gamma(2)$ and $\tilde\gamma(t)\neq \tilde\gamma(t+1)$  for every $t\in (0,1]$. This map induces a continuous closed curve with no self-intersections
\begin{eqnarray*}
\gamma:[0,2]&\to&\T\times \C\\
\gamma(\theta)&\mapsto&\big(\langle \theta\rangle ,\gamma(\theta)\big).
\end{eqnarray*} 
where $\langle \theta\rangle $ means the fractional part of $\theta$. In other words, $\gamma$ is a curve which turns twice in the $\T$ direction. We can think $\gamma$ as a $2-$fold covering of $\T$.
 We call such a function a $2-$curve. For example, $\gamma:[0,2]\to \C$ such that $\gamma(\theta)=ie^{\pi i \theta}$ is a $2$-curve.
\\
 
  For a    $2-$curve $\gamma$  we denote by $\gamma_{\theta}$ the fiber $\gamma\cap \big(\{ \theta\}\times \C\big)$. The choice of a point $(\theta,z)$ in a fiber $\gamma_{\theta}$ defines an order on $\gamma_{\theta}=\{z_{1,\theta}=(\theta,z),z_{2,\theta}\}$, induced by the positive direction on the circle. Moreover, for any other fiber $\gamma_{\tilde{\theta}}$ this choice induces an order $\gamma_{\tilde{\theta}}=\{z_{1,\tilde{\theta}},z_{2,\tilde{\theta}}\}$ where $z_{1,\tilde{\theta}}$ is the first point in $\gamma_{\tilde{\theta}}$ from $z_{1,\theta}$ following the positive orientation of the circle. \\
  
  A $2$-curve $\gamma$ is called {\it invariant} by $P$ if $P(\gamma)=\gamma$. In such case there exists $\tau\in \Z_2$ such that
    \[
  P(z_{j,\theta})=z_{j+\tau,\theta+\alpha}
  \]
   for every $\theta \in \T$ and every point $z_{j,\theta}$ in the fiber $\gamma_{\theta}$ (where first subscripts are taken in $\Z_2$). Note that an invariant $2$-curve have no counterpart in the non-fibred polynomial dynamics world. Indeed, the induced dynamics on $\gamma$ by every iterate $P^n$ is always minimal.    
   The real number
   \[
   \kappa_2(\gamma)=\frac{1}{2}\int_{\T}\log \big|\partial_zP_{\theta}(z_{1, \theta})\partial_zP_{\theta}(z_{2, \theta})\big|d\theta
   \]
is called the {\it multiplicator} of the invariant $2$-curve.
\\

   We define the {\it 2-unfolding} of $\gamma$ as the fibred quadratic polynomial
   \[
   _2P(t, z)=\left(t+\frac{\alpha+\tau}{2}, P_{2t}(z)\right).
   \]
   This dynamics consists just into two copies of $P$ and such that every half of the circle $\{t\in \T\}$ correspond to one circle $\{\theta\in \T\}$. The unfolding of the curve $\gamma$ gives raises to two disjoint invariant curves (turning just once in the $\T$ direction) for $_2P$:
  \begin{displaymath}
\gamma_1(t)=\left\{ \begin{array}{ccc}
  z_{1,2t}&\textrm{for }& 0\leq t\leq \frac{1}{2} \\
 z_{2, 2t} &\textrm{for }& \frac{1}{2}\leq t\leq 1,
 				\end{array}\right.
 \end{displaymath} 
and $\gamma_2(t)=\gamma_1(t+\frac{1}{2})$. The multiplicator of these  invariant curves verifies
\[
\kappa(\gamma_1)=\kappa(\gamma_2)=\kappa_2(\gamma).
\]

As in the non-fibred case, the quadratic family $\{z^2+\cc(\theta)\}_{\cc}$ is universal among all the fibred quadratic polynomials, as shown by the following:
  \begin{lema}[see Sester \cite{SEST99}] \label{lemasester1}
  Let $\alpha\in \T$, and three continuous functions $A:\T\to\C\setminus\{0\}$, $B:\T\to\C$, $C:\T\to \C$, which are the coefficients of the fibred quadratic polynomial
  \[
  P(\theta,z)=\big(\theta+\alpha,A(\theta)z^2+B(\theta)z+C(\theta)\big).
  \]
  There exists a continuous fibred transformation, affine at each fiber,
   \[
     W(\theta,z)=\big(\theta,u_1(\theta)z+u_2(\theta)\big),
  \]
and a continuous parameter $\cc:\T\to\C$ such that
\[
W\circ P\circ W^{-1}(\theta,z)=\big(\theta+\alpha,z^2+\cc(\theta)\big)\quad_{\blacksquare}
\]
  \end{lema}
\section{The example} 
 Let's study periodic points of exact period  $2$ for the quadratic polynomial $q_c(z)=z^2+c$. These points are given by the zeroes of the polynomial $Q_2(z,c)=z^2+z+c+1$, and  thus, they are
 \begin{equation*}
 z_1(c)=-\frac{1}{2}+\sqrt{-\frac{3}{4}-c}\quad , \quad   z_2(c)=-\frac{1}{2}-\sqrt{-\frac{3}{4}-c}.
 \end{equation*}
  Let $\cc:\T\to\C$ be a small continuous loop around the parabolic parameter $c=-3/4$. By following the two periodic  points $z_1(\cc(\theta)), z_2(\cc(\theta))$ when $\theta$ runs from $0$ to $1$, we see that one tour to the loop of parameters give raises to a transposition of the periodic  points of $q_{\cc(0)}$.  In fact, these two point move with $\theta$ drawing a $2$-curve, let's call it $\gamma$. If the loop $\cc$ is small enough, then it need to spend some time in the interior of the $2$-limb 
  \[
  \{c\quad \big|\quad |c-(-1)|<1/4\}
  \]
  which correspond to the parameters $c$ such that the cycle of period $2$ is attracting.  By reparametrizing if necessary, we can assume that the loop $\cc$ spend a lot of time inside this limb in such a way that
  \[
  \int_{\T}\log \big|q_{\cc(\theta)}'(z_1(\theta))q_{\cc(\theta)}'(z_2(\theta))\big|d\theta<0.
  \]
  Hence, the $2$-curve $\gamma$ is invariant by $P(\theta, z)=(\theta, z^2+\cc(\theta))$ and the multiplicator $\kappa_2(\gamma)$ is negative.  Also note that $P(z_{j}(\theta))=z_{j+1}(\theta)$ and then $\tau=1$. For $\alpha\in \R$ we define 
  \[
  P^{\alpha}(\theta, z)=\left(\theta+\alpha, z^2+\cc(\theta)\right).
  \]
  The $2$-curve $\gamma$ is no longer invariant by $P^{\alpha}$. We want to correct the quadratic fibre of the map $P^{\alpha}$ in order to left $\gamma$ invariant. We put 
  \[
  b_{\theta}(z)=z+\left(z_1(\theta+\alpha)-z_1(\theta)\right)\frac{z-z_2(\theta)}{z_1(\theta)-z_2(\theta)}+\left(z_2(\theta+\alpha)-z_2(\theta)\right)\frac{z-z_1(\theta)}{z_2(\theta)-z_1(\theta)}.
  \]
  This is an affine map, which depends continuously on $\theta$ and sends $z_j(\theta)$ to $z_j(\theta+\alpha)$. We construct the corrected fibred quadratic polynomial
  \[
  \tilde{P^{\alpha}}(\theta,z)=\left(\theta+\alpha, b_{\theta}\left(z^2+\cc(\theta)\right)\right)
  \] 
   which left the $2$-curve $\gamma$ invariant.  If we choose $\alpha$ small enough (and provided that $\cc$ is Lipschitz), then $\partial_zb_{\theta}$ is close to $1$ and the multiplicator of $\gamma$ (as an invariant curve for $\tilde{P^{\alpha}}$) is still negative.  By considering the $2$-unfolding of $\gamma$ (see previous section), we get a fibred quadratic polynomial with two attracting invariant curves $\gamma_1, \gamma_2$.  By Lemma \ref{lemasester1} we can assume that this example has the form $$P(\theta, z)=\left(\theta+\frac{\alpha+1}{2}, z^2+\tilde{\cc}(\theta)\right).$$

The circle $\{(\theta, 0)\}_{\theta\in \T}$ contains the critical point of each fibre. In the non-fibred quadratic setting  the critical point must be attracted to the (unique) attracting cycle. In a general fibred quadratic dynamics (in the form $z^2+\cc(\theta)$), having an  attracting invariant curve $\eta$, the set $\Omega(\eta)$ of $\theta\in \T$ such that the critical point $(\theta, 0)$ is attracted by $\eta$ is open. Moreover,  we have the following result due to Sester:
\begin{prop}[see Sester \cite{SEST99}, Proposition 3.2]
The set $\Omega(\eta)$ is non-empty $\quad_{\blacksquare}$
\end{prop}   
Since $\gamma_1\cap \gamma_2=\emptyset$, by connexity the set $\T\setminus\left(\Omega(\gamma_1)\cup\Omega(\gamma_2)\right)$ is non-empty.  
\\

\noindent{{\it Concluding remarks. }}
Any fibred quadratic polynomial close enough to $P$ still have two attracting invariant curves and the same remark on the critical points holds.  Let's suppose that $\hat P$ is close to $P$ and hyperbolic (that is, the post critical set does not accumulate the fibred Julia set, see \cite{JONS99, SEST99} for definitions).  In that case, there should be a  critical point $(\tilde \theta, 0)$ going to infinity or accumulating around  a non-attracting object (since the circle cannot be covered by disjoint open sets) in the fibred Fatou set. 
\bibliographystyle{plain}
  \bibliography{ejemplo_doble_curva.bib}

 \end{document}